\documentclass[12pt]{amsart}
\usepackage{amssymb,amsmath,amsthm,amscd}
\usepackage[all]{xy}

\numberwithin{equation}{section}

\newtheoremstyle{fancy1}{10pt}{10pt}{\itshape}{12pt}{\textsc\bgroup}{.\egroup}{8pt}{
}
\newtheoremstyle{fancy2}{10pt}{10pt}{}{12pt}{\itshape}{.}{8pt}{ }

\theoremstyle{fancy1}

\newtheorem{lem}[equation]{Lemma}
\newtheorem{prop}[equation]{Proposition}
\newtheorem{thm}[equation]{Theorem}

\newtheorem*{main*}{Result}
\newtheorem*{conjecture*}{Conjecture}
\newtheorem*{cor*}{Corollary}

\newcommand{\gal}[2]{\operatorname{Gal}\left( #1 / #2 \right)}

\newtheorem*{def*}{Definition}

\newtheorem*{rem*}{Remark}

\newtheorem*{example*}{Example}

\newtheorem*{examples*}{Examples}

\theoremstyle{remark}
\newtheorem*{case*}{Case}
\newtheorem*{case}{Case}

\newcommand{\cref}[1]{Corollary~\ref{#1}}

\newcommand{\Q}{{\mathbb{Q}}}

\newcommand{\Z}{{\mathbb{Z}}}
\def\con#1=#2(#3){#1 \equiv #2 \bmod{#3}}

\newcommand{\Aut}{\ensuremath{\operatorname{Aut}}}

\newcommand{\GL}{\ensuremath{\operatorname{GL}}}

\begin{document}

\title[Infinite class field towers]{Infinite class field towers}

\author{Jing Long Hoelscher}
\address{University of Arizona\\
     Tucson, AZ 85721}
\email{jlong@math.arizona.edu}

\begin{abstract}
This paper studies infinite class field towers of number fields $K$ that are ramified over $\Q$ only at one finite prime. In particular, we show the existence of such towers for a general family of primes including $p=2$, $3$ and $5$.

\end{abstract}

\maketitle

%-------------- Article Text--------------------

\section{Introduction}
This paper is concerned with number fields ramified only at one finite prime over $\Q$, which allow an infinite class tower. Suppose $K=K_0$ is an algebraic number field. For $i=0,1,2,\dots$, take $K_{i+1}$ to be the maximal abelian unramified extension of $K_i$, i.e. $K_{i+1}$ is the Hilbert class field of $K_i$. If there are no integers $i\geq 0$ such that $K_i=K_{i+1}$, we say that $K$ \emph{admits an infinite class field tower}. I. R. \v{S}afarevi\v{c} gave an example of an infinite class tower of a number field ramified at seven finite primes over $\Q$ in \cite{Sha}. More examples can be found in \cite{GS}, and \cite{Ha}, etc. Later Ren\'{e} Schoof showed in \cite{Sc} that $K=\Q(\zeta_{877})$, $\Q(\sqrt{-3321607})$, $\Q(\sqrt{39345017})$ and $\Q(\sqrt{-222637549223})$ each have an infinite class field tower where the only finite primes ramified in $K/\Q$ are respectively $877$, $3321607$, $39345017$, $222637549223$. In this paper we prove that either the cyclotomic field $K=\Q(\zeta_{p^m})$ or a degree $p$ extension $L$ of an unramified abelian extension $H$ of $K$ has an infinite class field tower under certain explicit conditions on the order of the class group of $K$. This yields a number field ramified only at one finite prime that has an infinite class field tower. The main theorem of this paper is:

\begin{thm}\label{maintheorem}
Let $p$ be a prime, suppose there exist a cyclotomic field $K=\Q(\zeta_{p^m})$ for some $m\in\mathbb{N}$ and a cyclic unramified Galois extension $H/K$ of prime degree $h$, which satisfy either of the following two conditions:
\begin{enumerate}
\item[(I)]$p\mid |Cl_H|$, $p$ is regular and $f_{p,h}^2-4f_{p,h}\geq 2h\cdot\varphi(p^m)$;
\item[(II)]$p\nmid |Cl_H|$ and $h\geq 2\varphi(p^{m+1})+4$,
\end{enumerate}
where $f_{p,h}$ is the order of $p$ in $(\Z/h\Z)^*$. Then there is a number field ramified only at $p$ and $\infty$ over $\Q$ that admits an infinite class tower.
\end{thm}
As an application we will show:
\begin{thm}\label{infiniteclassfield2}
For $p=2$, $3$ and $5$, there exist algebraic number fields ramified over $\Q$ only at $p$ and $\infty$ which admit an infinite class tower.
\end{thm}

\section{Infinite class field towers}
Let $K=\Q(\zeta_{p^m})$, and let $H$ be a cyclic unramified Galois extension over $K$ of prime degree $h$. Denote by $Cl_H$ the ideal class group of $H$. The principal prime ideal $(1-\zeta_{p^m})$ in $K$ above $p$ splits completely as $(1-\zeta_{p^m})=\prod_{i=1}^{h} \mathfrak{p}_i $ in the subfield $H$ of the Hilbert class field of $\Q(\zeta_{p^m})$. 

\begin{prop}\label{pdividclassnumber}
Under condition (I) in Theorem \ref{maintheorem}, the number field $H$ has an infinite class tower.
\end{prop}
\begin{proof}
Since $p||Cl_H|$, we can pick a degree $p$ cyclic unramified extension $M_0/H$. 
If $M_0/K$ is Galois, the Galois group $\gal{M_0}{K}\cong\Z/p\Z \rtimes \Z/h\Z$ since $(p,h)=1$. Now $p$ is regular, i.e. $p\nmid |Cl_{\Q(\zeta_p)}|$, by Theorem $10.4(a)$ in \cite{Wa} we know $p$ does not divide the class number of  $K=\Q(\zeta_{p^m})$. So the semi-direct product $Z/p\Z\rtimes \Z/h\Z$ is not abelian, we have a non-trivial homomorphism $\Z/h\Z\rightarrow \Z/(p-1)\Z$. So $h|(p-1)$, since $h$ is assumed to be a prime, i.e. the order $f_{(p,h)}$ of $p$ in $(Z/h\Z)^*$ is $1$, in contradiction to the condition $f_{(p,h)}^2-4f_{(p,h)}\geq 2h\varphi(p^m)$.

Otherwise $M_0/K$ is non-Galois. Take the Galois closure $\bar{M}_0$ of $M_0/K$, which is the composite of all conjugates of $M_0$ over $K$. The Galois group $\gal{\bar{M}_0}{H}$ is of the form $(\Z/p\Z)^l$ for some $l\in\mathbb{N}$, and $\gal{\bar{M}_0}{K}\cong (\Z/p\Z)^l\rtimes\Z/h\Z$. The semi-direct product gives a non-trivial homomorphism $\Z/h\Z\rightarrow \GL_l(\mathbb{F}_p)\cong \Aut((\Z/p\Z)^l)$, so $h\,|\,\hspace{0.1em}|\GL_l(\mathbb{F}_p)|=(p^l-1)(p^l-p)\cdots(p^l-p^{l-1})$. Then $h\,|\,(p^i-1)$, for some $i\leq l$, i.e.\, $p^i=1$ mod $h$, and the order $f_{(p,h)}$ of $p$ in $(\Z/h\Z)^*$ divides $i$. Thus $l\geq i\geq f_{p,h}$. So by the assumption we know $l^2-4l\geq 2h\cdot \varphi(p^m)$. On the other hand, there are $l$ linearly disjoint cyclic degree $p$ unramified extensions of $H$, i.e. $h_1=\dim H^1(G,\Z/p\Z)\geq l$, where $G=\gal{\Omega_H}{H}$ is the Galois group of the maximal unramified $p$-extension $\Omega_H$ of $H$. If $\Omega_H$ were finite, we would have $\frac{l^2}{4}-l\leq \frac{h_1^2}{4}-h_1 < h_2-h_1\leq r_1+r_2=h\varphi(p^m)/2$, where the first inequality comes from the fact $M_0/K$ is non-Galois thus $h_1\geq l\geq 2$. So $l^2-4l < 2h\varphi(p^m)$, contradiction.  
\end{proof}

From now on we will focus on the case under condition (II) in Theorem \ref{maintheorem}.
\begin{prop}\label{rayclass}
Assume $p\nmid |Cl_H|$, then there exist a prime $\mathfrak{p}$ in $H$ above $p$ and an integer $k$ such that $p$ divides $|Cl_H^{\mathfrak{p}^k}|/|Cl_H|$, where $Cl_H^{\mathfrak{p}^k}$ is the ray class group of $H$ for the modulus $\mathfrak{p}^k$.
\end{prop}
\begin{proof}Suppose the conclusion were false. Then for any integer $k$ and any prime ideal $\mathfrak{p}_i$ in $H$ above $p$, where $i=1,\cdots, h$, we have $p\nmid |Cl^{\mathfrak{p}_i^k}|/|Cl_H|$. Taking the modulus $\mathfrak{m}=\mathfrak{p}_1^k\cdots \mathfrak{p}_h^k$ for $k\in\mathbb{N}$, we have the following exact sequence from class field theory
$$1\longrightarrow \mathcal{O}^*/\mathcal{O}^{\mathfrak{m}}\longrightarrow (\mathcal{O}_H/\mathfrak{m})^*\longrightarrow Cl_H^{\mathfrak{m}}\longrightarrow Cl_H\longrightarrow 1,$$
where $\mathcal{O}^*$, resp.\ $\mathcal{O}^{\mathfrak{m}}$, is the group of units in $\mathcal{O}_H$, resp.\ of units $\equiv 1$ mod $\mathfrak{m}$ in $\mathcal{O}_H$. Combining with Chinese Remainder Theorem, we have 
\begin{equation}\label{rayclassgroup}
\frac{|Cl_H^{\mathfrak{m}}|}{|Cl_H|}=\frac{|(\mathcal{O}_H/\mathfrak{m})^*|}{|\mathcal{O}^*/\mathcal{O}^{\mathfrak{m}}|}=\frac{\prod_{i=1}^h|(\mathcal{O}_H/\mathfrak{p}_i^k)^*|}{\prod_{i=1}^h|\mathcal{O}^*/\mathcal{O}^{\mathfrak{p}_i^k}|}=\prod_{i=1}^h \frac{|(\mathcal{O}_H/\mathfrak{p}_i^k)^*|}{|\mathcal{O}^*/\mathcal{O}^{\mathfrak{p}_i^k}|}. 
\end{equation}
So $p\nmid \frac{|Cl_H^{\mathfrak{m}}|}{|Cl_H|}$. Taking $k=\varphi(p^m)k_0$ with $k_0\in\mathbb{N}$, the modulus $\mathfrak{m}=\mathfrak{p}_1^k\cdots\mathfrak{p}_h^k=(1-\zeta_{p^m})^k=(p)^{k_0}$ is a power of $(p)$. With the assumption $p\nmid |Cl_H|$, Equation \ref{rayclassgroup} says that there are no totally ramified abelian $p$-extensions of $H$ that are ramified over $\Q$ only at $p$. That is contrary to the fact that $H(\zeta_{p^a})/H$, with $a\in\mathbb{N}$ big enough, is a totally ramified $p$-extension ramified only at primes above $p$.  
\end{proof}

By Proposition \ref{rayclass} above, we know there exists an integer $k$ and a prime, say $\mathfrak{p}_1$ above $p$ in $H$, such that $p\mid \frac{|Cl_H^{\mathfrak{p}_i^k}|}{|Cl_H|}$. Since $p$ does not divide the class number $|Cl_H|$, there exists a degree $p$ extension $H_1$ of $H$ ramified only at $\mathfrak{p}_1$. Since $H_1/H$ is a Kummer extension, $H_1$ is of the form $H(\sqrt[p]{x_1})$ for some $x_1\in \mathcal{O}_H$. The conjugates $H_i$ of $H_1$ over $\Q(\zeta_{p^m})$ are degree $p$ extensions of $H$ ramified only at $\mathfrak{p}_i$, and they are of the form $H_i=H(\sqrt[p]{x_i})$ for some $x_i\in \mathcal{O}_H$. If $x_1\in \mathfrak{p}_1$, then $x_i\in \mathfrak{p}_i$ for $1\leq i\leq h$ since $x_i$ are the conjugates of $x_1$; similarly if $x_1\in \mathcal{O}^*_H$, we will have $x_i\in\mathcal{O}^*_H$ for any $1\leq i\leq h$. We denote by  $\delta_{H_{\mathfrak{p}_i}(\sqrt[p]{x_i})/H_{\mathfrak{p}_i}}$ and $\delta_{H_{\mathfrak{p}_i}(\sqrt[p]{x_1\cdots x_i})/H_{\mathfrak{p}_i}}$ respectively the discriminants of the $\mathfrak{p}_i$-adic extension $H_{\mathfrak{p}_i}(\sqrt[p]{x_i})/H_{\mathfrak{p}_i}$ and $H_{\mathfrak{p}_i}(\sqrt[p]{x_1\cdots x_h})/H_{\mathfrak{p}_i}$. Let $L$ be the field $H(\sqrt[p]{x_1\cdots x_h})$.

\begin{prop}\label{discriminant}
$\delta_{L/H}=\prod_{i=1}^h \delta_{H(\sqrt[p]{x_i})/H}.$
\end{prop}
\begin{proof}
The extension $H_i/H$ is ramified only at $\mathfrak{p}_i$, so the discriminant $\delta_{H_i/H}$ satisfies $\delta_{H_i/H}=\delta_{H_{\mathfrak{p}_i}(\sqrt[p]{x_i})/H_{\mathfrak{p}_i}}$. And for the discriminant of the extension $L/H$, we have $\delta_{L/H}=\prod_{i=1}^h \delta_{H_{\mathfrak{p}_i}(\sqrt[p]{x_1\cdots x_i})/H_{\mathfrak{p}_i}}.$ If suffices to show $\delta_{H_{\mathfrak{p}_i}(\sqrt[p]{x_1\cdots x_h})/H_{\mathfrak{p}_i}}=\delta_{H_{\mathfrak{p}_i}(\sqrt[p]{x_i})/H_{\mathfrak{p}_i}}.$ 

We may assume the valuation $\upsilon_i=\upsilon_{\mathfrak{p}_i\mathcal{O}_{H_{\mathfrak{p}_i}}}(x_i)\in\{0,1\}$ by Lemma $2.1$ of \cite{Da}. The extension $H_{\mathfrak{p}_i}(\sqrt[p]{x_j})/H_{\mathfrak{p}_i}$ is unramified at $\mathfrak{p}_i$ when $j\neq i$, by Theorem $2.4$ of \cite{Da} $x_j$ is a $p$ power of a unit mod $(\mathfrak{p}_i\mathcal{O}_{H_{\mathfrak{p}_i}})^{p^m}$. 

\begin{case}[a]For all $1\leq i\leq h$, $x_i\in\mathfrak{p}_i$. We have $x_1 \cdots x_h\in\mathfrak{p}_i$ and $\upsilon_{\mathfrak{p}_i\mathcal{O}_{H_{\mathfrak{p}_i}}}(x_i)=\upsilon_{\mathfrak{p}_i\mathcal{O}_{H_{\mathfrak{p}_i}}}(x_1\cdots x_n)$, since $x_j$ is a unit in $\mathcal{O}_{H_{\mathfrak{p}_i}}$ for $j\neq i$. By Theorem $2.4$ of \cite{Da} we know $\delta_{H_{\mathfrak{P}_i}(\sqrt[p]{x_1\cdots x_h})/H_{\mathfrak{p}_i}}=\delta_{H_{\mathfrak{P}_i}(\sqrt[p]{x_i})/H_{\mathfrak{p}_i}}.$ 
\end{case}
\begin{case}[b]For all $1\leq i\leq h$, $x_i\in \mathcal{O}^*_H$. We know $x_j$ is a $p$th power mod $(\mathfrak{p}_i\mathcal{O}_{H_{\mathfrak{p}_i}})^{p^m}$ for $j\neq i$. If we define $$ \kappa(x) = \max_{0\leq l \leq p^m} \{ l\,|\,\exists \gamma\in \mathcal{O}_{H_{\mathfrak{p}_i}}:\gamma^p=x \hspace{0.5em} \text{mod}\hspace{0.5em} (\mathfrak{p}_i\mathcal{O}_{H_{\mathfrak{p}_i}})^l\}$$ as a function of $x$, we have $\kappa(x_i)=\kappa(x_1\cdots x_h)$. By Theorem $2.4$ of \cite{Da}, again we get $\delta_{H_{\mathfrak{P}_i}(\sqrt[p]{x_1\cdots x_h})/H_{\mathfrak{p}_i}}=\delta_{H_{\mathfrak{P}_i}(\sqrt[p]{x_i})/H_{\mathfrak{p}_i}}.$
\end{case}
\end{proof}
We will show next that $L=H(\sqrt[p]{x_1\cdots x_h})$ has many linearly disjoint unramified abelian extensions.
\begin{lem}\label{unramified}
The extension $L(\sqrt[p]{x_i})/L$ is unramified for $1\leq i\leq h$.
\end{lem}
\begin{proof}Let $M=H(\sqrt[p]{x_1},\cdots,\sqrt[p]{x_h})$. We will show $M/L$ is unramified. Consider the discriminant of $M/H$. On the one hand, for each $1\leq i \le h$, $H_i/H$ is ramified only at $\mathfrak{p}_i$, so these extensions are linearly disjoint. The extension $M/H$, being the compositum of these extensions, has discriminant $\delta_{M/H}=\prod_{i=1}^h \delta_{H_i/H}^{p^{h-1}}$. On the other hand, considering the tower of extensions $M/L/H$, we have by Proposition \ref{discriminant}
$$\delta_{M/H}=\delta_{L/H}^{p^{h-1}}N_{L/H}(\delta_{M/L})=(\prod_{i=1}^h \delta_{H_i/H})^{p^{h-1}}N_{L/H}(\delta_{M/L}).$$
Comparing these two equations, we get $\delta_{M/L}=1$, i.e. $M/L$ is unramified.

\end{proof}
\begin{prop}\label{pnotdividclassnumber}
With the assumption of condition (II) in Theorem \ref{maintheorem}, the number field $L$ admits an infinite class tower.
\end{prop} 
\begin{proof}
Let $\Omega$ be the maximal unramified pro-$p$ extension of $L$, with Galois group $G=\gal{\Omega}{L}$. We claim $\Omega/L$ is infinite. Suppose $\Omega/L$ were a finite extension. Then $G$ is a finite $p$-group and $\Omega$ admits no cyclic unramified extensions of degree $p$. By Proposition $29$ in \cite{Sh} we know that 
\begin{equation}
h_2-h_1\leq r_1+r_2=r_2=\frac{p\varphi(p^m)h}{2}=\frac{\varphi(p^{m+1})h}{2},
\end{equation}
where $h_i=\dim H^i(G,\Z/p\Z)$, and $r_1,r_2$ denote the number of real places and pairs of complex places of $L$. By Lemma \ref{unramified}, we have $h$ linearly disjoint unramified $p$-extensions $L(\sqrt[p]{x_i})$ over $L$, so
\begin{equation}
h_1(G)\geq h.
\end{equation}
By the \v{S}afarevi\v{c}-Golod theorem (page $82$ of \cite{Sh}), we have
\begin{equation}
h_2>\frac{h_1^2}{4}.
\end{equation}
Condition (II) says $h$ is a prime with $h(\geq 2\varphi(p^{m+1})+4)>4$. Combining all the inequalities and the fact that the function $x^2-4x$ monotonically increases with $x$ if $x\geq 2$, we have $\frac{h^2}{4}-h\leq \frac{h_1^2}{4}-h_1<h_2-h_1\leq\frac{\varphi(p^{m+1})h}{2}$, which implies $h<2\varphi(p^{m+1})+4$, contradiction.
\end{proof}
Combining Proposition \ref{pnotdividclassnumber} and Proposition \ref{pdividclassnumber}, we can conclude Theorem \ref{maintheorem} in the introduction.

\section{Examples for small primes}
As a consequence of Theorem \ref{maintheorem}, we will verify in the cases $p=2$, $3$, $5$ that there are number fields ramified only at $p$ and $\infty$ which have an infinite class tower. The following examples rely on Table $3$ of \cite{Wa} for the relative class number of the cyclotomic number field $K$, i.e.\, the ratio of the class number $|Cl_K|$ of $K$ and the class number $|Cl_{K^+}|$ of the maximal real subfield $K^+$. Also the computation of the order $f_{(p,h)}$ of $p$ in $(\Z/h\Z)^*$ is done by \cite{PARI2}.
{\setlength{\tabcolsep}{0.1cm}
\renewcommand{\arraystretch}{1.6}
\begin{table}[!h]
\begin{center}
\begin{tabular}{|c||c|c|c|}
\hline
   & $K$ & $h$  & $f_{p,h}$  \\
\hline \hline
$p=2$ & $\Q(\zeta_{128})$ & $21121$ & $10560$  \\
\hline
$p=3$ & $\Q(\zeta_{81})$ & $2593$ &  $648$  \\
\hline
$p=5$ & $\Q(\zeta_{125})$ & $20602801$ & $10301400$  \\
\hline
\end{tabular}
\end{center}
\end{table}} 

\begin{case}[$p=2$]
Pick $K=\Q(\zeta_{128})$. The relative class number of $K$ is $359057=17\cdot 21121$. Let $H$ be a subfield of the Hilbert class field of $K$ with $\gal{H}{K}\cong \Z/21121\Z$. When $2\mid |Cl_H|$, condition (I) in Theorem \ref{maintheorem} is satisfied since $2$ is regular and $f_{2,21121}=10560$ thus $f_{2,21121}^2-4f_{2,21121}\geq 2\cdot h \cdot\varphi(128)$; when $2\nmid |Cl_H|$, condition (II) in Theorem \ref{maintheorem} is satisfied since $h=21121$ and thus $h\geq 2\varphi(2^8)+4=260$. 
\end{case}

\begin{case}[$p=3$]
Pick $K=\Q(\zeta_{81})$. The relative class number of $K$ is $2593$. Let $H$ an abelian unramified extension over $K$ with Galois group $\Z/2593\Z$. When $3\mid |Cl_H|$, condition (I) in Theorem \ref{maintheorem} is satisfied since $3$ is regular and $f_{3,2593}=648$ thus $f_{3,2593}^2-4f_{3,2593}\geq 2\cdot h\cdot \varphi(81)$; when $3\nmid |Cl_H|$, condition (II) in Theorem \ref{maintheorem} is satisfied since $h=2593$ thus $h\geq 2\varphi(3^4)+4=112$.
\end{case}

\begin{case}[$p=5$]
Pick $K=\Q(\zeta_{125})$. The relative class number of $K$ is $2801\cdot 20602801$. Let $H$ be the subfield of the Hilbert class field of $K$ with $\gal{H}{K}\cong \Z/20602801\Z$. When $5\mid |Cl_H|$, condition (I) in Theorem \ref{maintheorem} is satisfied since $5$ is regular and $f_{5,20602801}=103011400$ thus $f_{5,20602801}^2-4f_{5,20602801}\geq 2\cdot h \cdot\varphi(125)$; when $5\nmid |Cl_H|$, condition (II) in Theroem \ref{maintheorem} is satisfied since $h=20602801$ thus $h\geq 2\varphi(5^4)+4=1004$.
\end{case}
All the cases above complete the proof of Theorem \ref{infiniteclassfield2}.

\bigskip

\providecommand{\bysame}{\leavevmode\hbox
 to3em{\hrulefill}\thinspace}

 \end{document}